# CLARK–OCONE FORMULA AND VARIATIONAL REPRESENTATION FOR POISSON FUNCTIONALS

By Xicheng Zhang[1]

*University of New South Wales and
Huazhong University of Science and Technology*

In this paper we first prove a Clark–Ocone formula for any bounded measurable functional on Poisson space. Then using this formula, under some conditions on the intensity measure of Poisson random measure, we prove a variational representation formula for the Laplace transform of bounded Poisson functionals, which has been conjectured by Dupuis and Ellis [*A Weak Convergence Approach to the Theory of Large Deviations* (1997) Wiley], p. 122.

**1. Introduction.** Let $W$ be a standard $d$-dimensional Brownian motion. The following elegant formula for the Laplace transform of a bounded and measurable functional $F$ of Brownian motion was first established by Boué and Dupuis [1]:

$$(1) \qquad -\log \mathbb{E}[e^{-F}] = \inf_v \mathbb{E}\bigg[F\bigg(\cdot + \int_0^\cdot v_s\,ds\bigg) + \frac{1}{2}\int_0^1 |v_s|^2\,ds\bigg],$$

where the infimum is taken for all processes $v$ that are progressively measurable with respect to the augmented filtration generated by Brownian motion. This result was later extended to Hilbert space-valued Brownian motion by Budhiraja and Dupuis [3]. Furthermore, the author in [21] extended this representation to the abstract Wiener space, and gave a simplified proof by using Clark–Ocone's formula. This formula has proven to be useful in deriving various asymptotic results in large deviations (cf. [1, 2, 3, 16, 17, 18]).

For Poisson functionals, a similar representation formula has been conjectured by Dupuis and Ellis in [4], page 122, from the background of control

Received February 2008.
[1]Supported by ARC Discovery Grant DP0663153 of Australia.

theory:

$$-\log \mathbb{E}[e^{-F}]$$
$$= \inf_\phi \mathbb{E}\bigg[\int_0^1 \int_{\mathbb{R}^d} [\phi(y,t)\log\phi(y,t) - \phi(y,t) + 1]\nu_{\bar{X}(t)}(dy)\,dt + F(\bar{X})\bigg],$$

where the infimum is taken over all suitable controls $\phi$, and $\bar{X}$ is a controlled Markov process with jump defined by the generator

$$\int_{\mathbb{R}^d} [f(x+y) - f(x)]\phi(y,t)\nu_x(dy).$$

Here, $\nu_x(dy)$ is the jump intensity of a Markov process.

However, there is no rigorous proof for this variational formula up to now. In the present paper we will attempt to give a rigorous proof in a more general setting. Roughly speaking, let $(\Omega, \mathbb{P})$ be the canonical Poisson space (simple configuration space over $[0,1]\times\mathbb{R}^d$) and $\nu$ an intensity measure on $\mathbb{R}^d$. For any bounded random variable $F$ on $\Omega$, we want to prove that

$$-\log\mathbb{E}(e^{-F}) = \inf_\phi \mathbb{E}\bigg[\int_0^1\int_{\mathbb{R}^d}[\phi(y,t)\log\phi(y,t) - \phi(y,t) + 1]\nu(dy)\,dt + F\circ\Gamma_\phi^-\bigg],$$

where the infimum runs over some classes of predictable processes, and $\Gamma_\phi^-$ is a predictable transformation on $\Omega$ associated with $\phi$.

In contrast to the Wiener space case, the main difficulty of proving this formula comes from the nonlinearity of Poisson space. In particular, the Girsanov theorem for the Poisson measure is related to some nonlinear invertible and predictable transformations on $\mathbb{R}^d$ (cf. [5], Theorem 3.10.21). Indeed, the definition of the above $\Gamma_\phi^-$ depends on solving a mass transportation problem or the classical Monge–Ampère equation. More precisely, to a given positive function $\phi$, we need to seek an invertible transformation $x \mapsto y(x)$ of $\mathbb{R}^d$ such that, for all test functions $f \in C_0(\mathbb{R}^d)$,

$$\int_{\mathbb{R}^d} f(y(x))\nu(dx) = \int_{\mathbb{R}^d} f(x)\phi(x)\nu(dx),$$

which is formally equivalent to solving the following nonlinear PDE in the case of $\nu(dx) = \theta(x)\,dx$:

$$\theta(y^{-1}(x)) \cdot \det(\nabla y^{-1}(x)) = \theta(x)\phi(x).$$

For an optimal mass transportation problem, we refer to the book of Villani [20]. Since our problem has no constraint conditions on $y$, an easy solution can be constructed when $\nu$ has full support and no charges on $d-1$-dimensional subspaces, and satisfies an extra mild assumption. More detailed discussions are given in Section 5.

In order to prove the above variational representation formula, the first step is to establish the following Clark–Ocone formula: for any bounded functional $F$,

$$F = \mathbb{E}F + \int_0^1 \int_{\mathbb{U}} {}^p D_{(u,t)} F \tilde{\mu}(du, dt),$$

where $\tilde{\mu}$ is the compensated Poisson random measure, ${}^p D_{(u,t)} F$ is the predictable projection of $D_{(u,t)} F$ and $D$ is the difference operator [see (6) below]. The proof of this formula depends on an integration by parts formula given in Picard [11, 12]. Although there are many martingale representation formulas for Poisson functionals (e.g., see [8, 9, 13]), the well-known results are mainly concentrated on the representation for functionals in the first order Sobolev space by using the Chaos decomposition. The main point for us is that ${}^p D_{(u,t)} F$ is a bounded predictable process.

This paper is organized as follows: In Section 2 some notation and necessary lemmas are given as preliminaries. In Section 3 we prove the Clark–Ocone formula for bounded Poisson functionals. In Section 4 we shall prove two variational representation formulas for Poission functionals. One (Theorem 4.4 below) is weaker, and needs no assumption. Another (Theorem 4.11 below) is stronger, and needs to work in a locally compact metric space, and also requires some extra assumptions [see (H1) and (H2) below]. In Section 5 we discuss these two extra assumptions, and give a solution when $\mathbb{U} = \mathbb{R}^d$ and the intensity $\nu$ satisfies certain assumptions.

**2. Preliminaries.** Let $\mathbb{U}$ be a Lusin space, that is, a Hausdorff space that is the image of a Polish space under a continuous bijection. We fix a $\sigma$-finite and infinite measure $\nu$ on $(\mathbb{U}, \mathcal{B}(\mathbb{U}))$. Since $\mathbb{U} \times [0,1]$ is still a Lusin space and has the same cardinality with $\mathbb{R}$, it is well known that $(\mathbb{U} \times [0,1], \mathcal{B}(\mathbb{U} \times [0,1]))$ is isomorphic to $([0,1], \mathcal{B}([0,1]))$ (cf. [6], Proposition 8.6.12, and [10], Theorem 2.12). This property was used in the proof of [12], Theorem 1, and so in [11], Lemma 1.4 (see Theorem 3.2 below), which will be our basis for subsequent proofs.

Let $\Omega$ be the space of all integer-valued measures $\omega$ on $\mathbb{U} \times [0,1]$ such that $\omega(\{u,t\}) \leq 1$ for any $(u,t) \in \mathbb{U} \times [0,1]$, and $\omega(A \times [0,1]) < +\infty$ for any $A \in \mathcal{B}(\mathbb{U})$ with $\nu(A) < \infty$. The canonical random measure on $\Omega$ is then defined by

$$\mu_\omega(A \times (0,t]) := \omega(A \times (0,t]), \qquad t \in [0,1], A \in \mathcal{B}(\mathbb{U}).$$

The filtration $(\mathcal{F}_t)_{t \in [0,1]}$ is defined by

$$\mathcal{F}_t := \sigma\{\mu_\omega(A \times (0,s]) : s \leq t, A \in \mathcal{B}(\mathbb{U})\}.$$

We shall simply write $\mathcal{F}_1$ as $\mathcal{F}$. Let $\mathbb{P}$ be the probability measure on $(\Omega, \mathcal{F})$ such that $\mu_\omega$ is a Poisson random measure with the intensity measure $\nu(du)$.

That is, for any $A \in \mathcal{B}(\mathbb{U})$ and $t \in [0,1]$, the random variable $\omega \mapsto \mu_\omega(A \times (0,t])$ is a Poisson random variable with mean $\nu(A) \cdot t$, and $\omega \mapsto \mu_\omega(I_i \times A_j)$ are independent if the sets $I_i \times A_j$ are disjoint. We shall also denote by $\tilde{\mu}_\omega$ the compensated Poisson random measure $\mu_\omega - \pi$, where $\pi(du, dt) := \nu(du) \times dt$, and $dt$ is the Lebesgue measure on $[0,1]$.

Let $\mathcal{F}_t^\mathbb{P}$ be the completion of $\mathcal{F}_t$ with respect to $\mathbb{P}$, then $(\Omega, \mathcal{F}^\mathbb{P}, \mathbb{P}; (\mathcal{F}_t^\mathbb{P})_{t \in [0,1]})$ forms a complete filtration probability space. We shall denote by $\mathcal{P}$ the predictable $\sigma$-field associated with $(\mathcal{F}_t^\mathbb{P})_{t \in [0,1]}$, which is generated by all left continuous $\mathcal{F}_t^\mathbb{P}$-adapted processes. For the simplicity of notation, we shall write for $p \in [1, \infty]$

$$\mathbb{L}^p := L^p(\mathbb{U} \times [0,1] \times \Omega, \mathcal{B}(\mathbb{U} \times [0,1]) \times \mathcal{F}^\mathbb{P}, \pi \times \mathbb{P})$$

and

$$\mathbb{L}_\mathcal{P}^p := L^p(\mathbb{U} \times [0,1] \times \Omega, \mathcal{B}(\mathbb{U}) \times \mathcal{P}, \nu \times dt \times \mathbb{P}).$$

Let $\mathscr{C}$ be the linear span of the following simple processes:

$$\phi(u, t, \omega) := 1_{(t_0, t_1]}(t) \cdot g(u, \omega), \tag{2}$$

where $0 \leq t_0 < t_1 \leq 1$ and $g$ is bounded and $\mathcal{B}(\mathbb{U}) \times \mathcal{F}_{t_0}$-measurable and satisfies

$$g(u, \omega) \cdot 1_{U^c}(u) = 0 \qquad \text{for some } U \in \mathcal{B}(\mathbb{U}) \text{ with } \nu(U) < +\infty. \tag{3}$$

REMARK 2.1. For $g \in \mathcal{B}(\mathbb{U}) \times \mathcal{F}_{t_0}^\mathbb{P}$, by the monotone class theorem, we can find a $\tilde{g} \in \mathcal{B}(\mathbb{U}) \times \mathcal{F}_{t_0}$ such that $\tilde{g} = g$, $\nu \times \mathbb{P}$-a.s.

The following lemma is standard. The construction will also be used in the proof of Lemma 4.8 below.

LEMMA 2.2. $\mathscr{C}$ is dense in $\mathbb{L}_\mathcal{P}^p$ for any $p \in [1, \infty)$.

PROOF. We sketch the proof. Let $\phi \in \mathbb{L}_\mathcal{P}^p$. For $\varepsilon \in (0, 1/2)$, we first extend $\phi$ to $[-2\varepsilon, 0]$ by setting $\phi(u, t, \omega) = 0$ for $t \in [-2\varepsilon, 0]$, and then define

$$\phi_\varepsilon(u, t, \omega) := \frac{1}{\varepsilon^2} \int_{t-\varepsilon}^t \int_{s-\varepsilon}^s \phi(u, r, \omega) \, dr \, ds, \qquad t \in [0,1].$$

Obviously, $t \mapsto \phi_\varepsilon(u, t, \omega)$ is a continuous differentiable and $\mathcal{F}_t^\mathbb{P}$-adapted process, and satisfies

$$\int_0^1 |\phi_\varepsilon(u, t, \omega)|^p \, dt \leq \int_0^1 |\phi(u, t, \omega)|^p \, dt,$$

$$\int_0^1 |\phi_\varepsilon'(u, t, \omega)|^p \, dt \leq \frac{2^{p+1}}{\varepsilon^p} \int_0^1 |\phi(u, t, \omega)|^p \, dt.$$

Second, for $\varepsilon \in (0, 1/2)$ and $n \in \mathbb{N}$, we define

$$\phi_{\varepsilon,n}(u, t, \omega) := \sum_{k=0}^{n-1} 1_{(kn^{-1}, (k+1)n^{-1}]}(t) \cdot \phi_\varepsilon(u, kn^{-1}, \omega).$$

Then

$$\int_0^1 |\phi_{\varepsilon,n}(u, t, \omega)|^p \, dt \le \sup_{t \in [0,1]} |\phi_\varepsilon(u, t, \omega)|^p \le \int_0^1 |\phi'_\varepsilon(u, t, \omega)|^p \, dt.$$

Last, let $(U_m)_{m \in \mathbb{N}}$ be an increasing sequence of Borel subsets of $\mathbb{U}$ such that $\bigcup_m U_m = \mathbb{U}$ and $\nu(U_m) < +\infty$, and define

$$\phi_\varepsilon^m(u, kn^{-1}, \omega) := (-m) \vee (\phi_\varepsilon(u, kn^{-1}, \omega) \wedge m) \cdot 1_{U_m}(u).$$

By the diagonalization method and the dominated convergence theorem, we may find the desired approximation in $\mathbb{L}_{\mathcal{P}}^p$ by Remark 2.1. □

We recall the notion about the relative entropy as follows (cf. [4]).

DEFINITION 2.3. Let $\mathscr{P}(\Omega)$ denote the set of all probability measures defined on $(\Omega, \mathcal{F})$. For $\gamma \in \mathscr{P}(\Omega)$, the relative entropy function $R(\cdot \| \gamma)$ is a mapping from $\mathscr{P}(\Omega)$ into $\mathbb{R} \cup \infty$ given by

$$R(\gamma' \| \gamma) := \mathbb{E}^{\gamma'}\left( \log \frac{d\gamma'}{d\gamma} \right),$$

whenever $\gamma' \in \mathscr{P}(\Omega)$ is absolutely continuous with respect to $\gamma$ such that the above integral is finite, where $\mathbb{E}^{\gamma'}$ denotes the expectation with respect to $\gamma'$. In all other cases, $R(\gamma' \| \gamma) := \infty$.

The following proposition can be found in [4], Proposition 1.4.2.

PROPOSITION 2.4. *Let $\gamma \in \mathscr{P}(\Omega)$, and $F$ a bounded random variable on $(\Omega, \mathcal{F})$.*

  (i) *We have the following variational formula:*

$$-\log \mathbb{E}^\gamma(e^{-F}) = \inf_{\gamma' \in \mathscr{P}(\Omega)} [R(\gamma' \| \gamma) + \mathbb{E}^{\gamma'}(F)].$$

  (ii) *The infimum in* (i) *is uniquely attained at the probability measure $\gamma_0$ defined by*

(4) $$\gamma_0(d\omega) = e^{-F(\omega)} / \mathbb{E}^\gamma(e^{-F}) \cdot \gamma(d\omega).$$

The following contents will be used in the second part of Section 4. In order to prove the variational representation formula for Poisson functionals, we need to endow $\Omega$ with a suitable topology such that $\Omega$ becomes a Polish space. For this, we assume that $\mathbb{U}$ is a noncompact locally compact connected complete metric space, and $\nu$ is a Radon measure on $\mathbb{U}$. Let $C_c(\mathbb{U} \times [0,1])$ denote the set of all continuous functions on $\mathbb{U} \times [0,1]$ with compact supports. The topology on $\Omega$ is taken as the weakest topology such that, for any $f \in C_c(\mathbb{U} \times [0,1])$, the function

$$(5) \qquad \omega \mapsto \langle f, \mu_\omega \rangle := \int_0^1 \int_\mathbb{U} f(u,t) \mu_\omega(du, dt) = \sum_{(u,t) \in \mathrm{supp}(\omega)} f(u,t)$$

is continuous, where $\mathrm{supp}(\omega)$ is the support of integer-valued Radon measure $\omega$, and the sum only has finite terms. By [15], Theorem 1.8, $\Omega$ is a Polish space under the above topology.

The following result can be found in [1], Lemma 2.8.

LEMMA 2.5. *Let $\gamma \in \mathscr{P}(\Omega)$ and $\{\gamma_n, n \in \mathbb{N}\} \subset \mathscr{P}(\Omega)$ satisfy*

$$\sup_{n \in \mathbb{N}} R(\gamma_n \| \gamma) < +\infty.$$

(i) *If $\{F_k, k \in \mathbb{N}\}$ is a sequence of uniformly bounded random variables converging to $F$, $\gamma$-a.s., then*

$$\lim_{k \to \infty} \sup_{n \in \mathbb{N}} \mathbb{E}^{\gamma_n} |F_k - F| = 0.$$

(ii) *If $\gamma_n$ converges weakly to the probability measure $\gamma$, then for any bounded random variable $F$ on $(\Omega, \mathcal{F})$*

$$\lim_{n \to \infty} \mathbb{E}^{\gamma_n}(F) = \mathbb{E}^\gamma(F).$$

Let $\mathcal{C}$ be the set of all cylindrical functions on $\Omega$ with the form

$$F(\omega) := h(\langle f_1, \mu_\omega \rangle, \ldots, \langle f_n, \mu_\omega \rangle), \qquad h \in C_c^\infty(\mathbb{R}^n), f_i \in C_c(\mathbb{U} \times [0,1]).$$

We also need the following standard result.

LEMMA 2.6. *Let $F$ be a bounded random variable on $(\Omega, \mathcal{F}^\mathbb{P}, \mathbb{P})$. Then there exists a family of functions $F_n \in \mathcal{C}$ with $\sup_n \|F_n\|_\infty \leq \|F\|_\infty$ such that for $\mathbb{P}$-almost all $\omega$*

$$F_n(\omega) \to F(\omega), \qquad as\ n \to \infty.$$

PROOF. We sketch it. Let $C_0(\mathbb{U} \times [0,1])$ be the completion of $C_c(\mathbb{U} \times [0,1])$ with respect to the uniform norm, which is then a separable Banach space. Let $\{f_i, i \in \mathbb{N}\} \subset C_c(\mathbb{U} \times [0,1])$ be a countable dense subset of $C_0(\mathbb{U} \times [0,1])$. We define $\mathcal{Q}_n := \sigma\{\langle f_i, \mu_\omega \rangle, i = 1, \ldots, n\}$. Then $\mathcal{Q}_n^{\mathbb{P}} \uparrow \mathcal{Q}_\infty^{\mathbb{P}} = \mathcal{F}^{\mathbb{P}}$. First, let $G_n := \mathbb{E}(F|\mathcal{Q}_n)$, then there exists some bounded measurable function $g_n$ on $\mathbb{R}^n$ such that

$$G_n(\omega) = h_n(\langle f_1, \mu_\omega \rangle, \ldots, \langle f_n, \mu_\omega \rangle).$$

Next, using the usual localizing and mollifying techniques, we may approach $h_n$ by $h_{n,k} \in C_c^\infty(\mathbb{R}^n)$. By the diagonalization method and extracting some subsequence if necessary, we then get the desired approximation sequence. □

**3. Clark–Ocone formula.** Let us first recall some definitions about the difference operator given in [11, 12]. For a fixed $(u,t) \in \mathbb{U} \times [0,1]$, define the transformation $\varepsilon_{(u,t)}^-$ and $\varepsilon_{(u,t)}^+$ on $\Omega$ by removing and adding a mass as follows: for $A \in \mathcal{B}(\mathbb{U} \times [0,1])$,

$$(\varepsilon_{(u,t)}^- \omega)(A) := \omega(A \cap \{(u,t)\}^c)$$

and

$$(\varepsilon_{(u,t)}^+ \omega)(A) := (\varepsilon_{(u,t)}^- \omega)(A) + 1_A(u,t).$$

It is clear that $(u,t,\omega) \mapsto \varepsilon_{(u,t)}^\pm \omega$ are $\mathcal{B}(\mathbb{U} \times [0,1]) \times \mathcal{F}^{\mathbb{P}}/\mathcal{F}^{\mathbb{P}}$-measurable.

For a functional $F$ on $\Omega$, the difference operator $D$ is defined by

(6) $$D_{(u,t)}F(\omega) := F \circ \varepsilon_{(u,t)}^+(\omega) - F(\omega).$$

Clearly, it is well defined except on a $\pi \times \mathbb{P}$-null set $N$. In the following, we always put $D_{(u,t)}F(\omega) = 0$ for $(u,t,\omega) \in N$. For a $\phi \in \mathbb{L}^1$, the divergence operator $\delta$ is defined by

$$\delta(\phi)(\omega) := \int_0^1 \int_{\mathbb{U}} \phi(u,t,\varepsilon_{(u,t)}^- \omega) \tilde{\mu}_\omega(du,dt).$$

We need the following simple lemma.

LEMMA 3.1. *For any $\phi \in \mathscr{C}$, we have*

$$\delta(\phi)(\omega) = \int_0^1 \int_{\mathbb{U}} \phi(u,t,\omega) \tilde{\mu}_\omega(du,dt).$$

PROOF. Let $\phi \in \mathscr{C}$ have the form (2). Notice that for any $A \in \mathcal{F}_{t_0}$ and $t > t_0, u \in \mathbb{U}$,

$$1_A \circ \varepsilon_{(u,t)}^- = 1_A.$$

Since $g$ is bounded and $\mathcal{B}(\mathbb{U}) \times \mathcal{F}_{t_0}$-measurable, by the monotone class theorem, we have for any $t > t_0$

$$g(u, \varepsilon^-_{(u,t)}\omega) = g(u, \omega).$$

Hence,

$$\begin{aligned}\delta(\phi)(\omega) &= \int_0^1 \int_{\mathbb{U}} 1_{(t_0,t_1]}(t) \cdot g(u, \varepsilon^-_{(u,t)}\omega)\tilde{\mu}_\omega(du, dt) \\ &= \int_0^1 \int_{\mathbb{U}} 1_{(t_0,t_1]}(t) \cdot g(u, \omega)\tilde{\mu}_\omega(du, dt) \\ &= \int_0^1 \int_{\mathbb{U}} \phi(u, t, \omega)\tilde{\mu}_\omega(du, dt).\end{aligned}$$

The result follows. $\square$

The following integration by parts formula can be found in [11], Lemma 1.4.

THEOREM 3.2. *Let $\phi \in \mathbb{L}^1$ and $F$ be a bounded random variable. Then*

$$\mathbb{E}(F\delta(\phi)) = \mathbb{E}\bigg(\int_0^1 \int_{\mathbb{U}} D_{(u,t)}F \cdot \phi(u,t)\pi(du,dt)\bigg).$$

Before proving the Clark–Ocone formula, we recall the following classical predictable projection theorem (cf. [19], page 173, Theorem 5.6).

LEMMA 3.3. *Let $\psi$ be a bounded measurable process on $\mathbb{U} \times [0,1] \times \Omega$. There exists a unique (up to indistinguishability with respect to $t$ for each $u$) predictable process $\phi \in \mathbb{L}^\infty_\mathcal{P}$ such that for every predictable stopping time $\tau$ and $u \in \mathbb{U}$*

(7) $\qquad \mathbb{E}(\psi(u,\tau) \cdot 1_{\{\tau < \infty\}} | \mathcal{F}^{\mathbb{P}}_{\tau-}) = \phi(u,\tau) \cdot 1_{\{\tau < \infty\}}, \qquad \mathbb{P}\text{-}a.s.$

*We shall write $\phi$ as $^p\psi$, which is called the predictable projection of $\psi$.*

PROOF. (Uniqueness) Let $\phi_1$ and $\phi_2$ be two predictable projections of $\psi$. Set

$$A := \{(u,t,\omega) : \phi_1(u,t,\omega) \neq \phi_2(u,t,\omega)\}.$$

Then for each $u \in \mathbb{U}$, the section $\Pi_u(A) := \{(t,\omega) : (u,t,\omega) \in A\} \in \mathcal{P}$. By the section theorem (cf. [19], page 172, Theorem 5.5) and (7), we have $\mathbb{P}(\Pi(\Pi_u(A))) = 0$, where $\Pi(\Pi_u(A)) = \{\omega : (t,\omega) \in \Pi_u(A), \exists t \in [0,1]\}$. Hence, for every $u \in \mathbb{U}$, $\phi_1(u,\cdot,\cdot)$ and $\phi_2(u,\cdot,\cdot)$ are indistinguishability.

(Existence) Let $\mathcal{M}$ be the class of all bounded measurable processes $\psi$ possessing a predictable projection. It is clear that $\mathcal{M}$ is a vector space containing the constants. Moreover, $\mathcal{M}$ is also a monotone class. In fact, let $\psi_n \in \mathcal{M}$ be a uniformly bounded increasing sequence with limit $\psi$. Let $\phi_n$ be the corresponding predictable projection of $\psi_n$. It is then easily checked by the monotone convergence theorem that $\overline{\lim}\,\phi_n$ is the predictable projection of $\psi$.

Hence, it is enough to prove that $\mathcal{M}$ contains all the processes of the form $\psi(u,t,\omega) = 1_{[0,t_0]}(t) \cdot g(u,\omega)$, which generates the $\sigma$-field $\mathcal{B}(\mathbb{U}) \times \mathcal{B}([0,1]) \times \mathcal{F}^{\mathbb{P}}$, where $g$ is bounded and $\mathcal{B}(\mathbb{U}) \times \mathcal{F}^{\mathbb{P}}$-measurable. Define
$$\phi(u,t) = 1_{[0,t_0]}(t) \cdot \mathbb{E}(g(u)|\mathcal{F}^{\mathbb{P}}_{t-}).$$
By Doob's optional stopping theorem, one then finds that such $\phi$ is a predictable projection of $\psi$. The proof is complete. $\square$

We now prove the following Clark–Ocone formula.

THEOREM 3.4. *Let $F$ be any bounded random variable on $\Omega$. Then*

(8) $$F = \mathbb{E}F + \int_0^1 \int_\mathbb{U} {}^pD_{(u,t)}F \tilde{\mu}(du,dt),$$

*where ${}^pD_{(u,t)}F \in \mathbb{L}^2_\mathcal{P} \cap \mathbb{L}^\infty_\mathcal{P}$ is the predictable projection of $D_{(u,t)}F$. Moreover,*
$$\mathbb{E}\left(\int_0^1 \int_\mathbb{U} |{}^pD_{(u,t)}F|^2 \pi(du,dt)\right)^2 < +\infty.$$

PROOF. It is well known that there exists a predictable process $\varphi \in \mathbb{L}^2_\mathcal{P}$ such that (cf. [7])
$$F = \mathbb{E}(F) + \int_0^1 \int_\mathbb{U} \varphi(u,t)\tilde{\mu}(du,dt).$$
By Lemma 3.1 and the isometry formula of the stochastic integral, we have for any $\phi \in \mathscr{C}$

(9) $$\mathbb{E}(F\delta(\phi)) = \mathbb{E}\left(\int_0^1 \int_\mathbb{U} \varphi(u,t) \cdot \phi(u,t)\pi(du,dt)\right).$$

On the other hand, by Theorem 3.2 and Fubini's theorem, we have for any $\phi \in \mathscr{C}$

(10)
$$\mathbb{E}(F\delta(\phi)) = \mathbb{E}\left(\int_0^1 \int_\mathbb{U} D_{(u,t)}F \cdot \phi(u,t)\pi(du,dt)\right)$$
$$= \int_0^1 \int_\mathbb{U} \mathbb{E}(\mathbb{E}(D_{(u,t)}F|\mathcal{F}^{\mathbb{P}}_{t-}) \cdot \phi(u,t))\pi(du,dt),$$
$$[\text{by } (7)] = \int_0^1 \int_\mathbb{U} \mathbb{E}({}^pD_{(u,t)}F \cdot \phi(u,t))\pi(du,dt)$$
$$= \mathbb{E}\left(\int_0^1 \int_\mathbb{U} {}^pD_{(u,t)}F \cdot \phi(u,t)\pi(du,dt)\right).$$

The formula (8) now follows by combining (9), (10) and Lemma 2.2.

By BDG's inequality (cf. [5], Theorem 4.1.12) and (8), we have

$$\mathbb{E}\left(\int_0^1\int_\mathbb{U}|{}^pD_{(u,t)}F|^2\mu(du,dt)\right)^2 \leq C\mathbb{E}\left(\int_0^1\int_\mathbb{U}{}^pD_{(u,t)}F\tilde{\mu}(du,dt)\right)^4$$
$$\leq C\mathbb{E}(F-\mathbb{E}F)^4,$$

where $C$ is a universal constant. Hence,

$$\mathbb{E}\left(\int_0^1\int_\mathbb{U}|{}^pD_{(u,t)}F|^2\pi(du,dt)\right)^2 \leq C\mathbb{E}\left(\int_0^1\int_\mathbb{U}|{}^pD_{(u,t)}F|^2\tilde{\mu}(du,dt)\right)^2$$
$$+ C\mathbb{E}\left(\int_0^1\int_\mathbb{U}|{}^pD_{(u,t)}F|^2\mu(du,dt)\right)^2$$
$$\leq C\mathbb{E}\left(\int_0^1\int_\mathbb{U}|{}^pD_{(u,t)}F|^4\mu(du,dt)\right)$$
$$+ C\mathbb{E}(F-\mathbb{E}F)^4$$
$$= C\mathbb{E}\left(\int_0^1\int_\mathbb{U}|{}^pD_{(u,t)}F|^4\pi(du,dt)\right)$$
$$+ C\mathbb{E}(F-\mathbb{E}F)^4 < +\infty.$$

The proof is complete. □

REMARK 3.5. In general, it is not known whether $\mathbb{E}(D_{(u,t)}F|\mathcal{F}_{t-}^\mathbb{P})$ is predictable, although for fixed $(u,t)\in\mathbb{U}\times[0,1]$, $\mathbb{E}(D_{(u,t)}F|\mathcal{F}_{t-}^\mathbb{P}) = {}^pD_{(u,t)}F$ a.s. by (7). However, Løkka in [8], Theorem 7 and Proposition 10, proved that for Lévy functional $F$, if $F$ belongs to the first-order Sobolev space, $(u,t)\mapsto\mathbb{E}(D_{(u,t)}F|\mathcal{F}_{t-}^\mathbb{P})$ is predictable. Compared with Løkka's result, Theorem 3.4 only requires that $F$ is bounded, and more importantly, the bound of ${}^pD_{(u,t)}F$ can be explicitly calculated from $F$, which is crucial for the next section. Moreover, this would also have some applications in mathematical finance as in [8], Section 5.

**4. Variational representation formula.** We begin with the following elementary lemma.

LEMMA 4.1. *Let $c_0 > -1$. Then for some $C > 0$ and any $x \geq c_0$,*

(11) $$|\log(1+x)| \leq C|x|, \qquad |\log(1+x) - x| \leq C|x|^2$$

*and*

(12) $$|(1+x)\log(1+x) - x| \leq C|x|^2.$$

Let $\phi \in \mathbb{L}^2_{\mathcal{P}}$ be a bounded predictable process satisfying

(13) $$\phi(u,t,\omega) \geq c_\phi > -1$$

and

(14) $$\mathbb{E}\left(\int_0^1 \int_{\mathbb{U}} \phi(u,t)^2 \pi(du,dt)\right)^2 < +\infty,$$

and such that $t \mapsto \mathscr{E}_t(\phi)$ is a square integrable $\mathcal{F}_t^{\mathbb{P}}$-martingale, where

(15) $$\begin{aligned}\mathscr{E}_t(\phi) := \exp\bigg\{&\int_0^t \int_{\mathbb{U}} \log(1+\phi(u,s))\tilde{\mu}(du,ds) \\ &+ \int_0^t \int_{\mathbb{U}} [\log(1+\phi(u,s)) - \phi(u,s)]\pi(du,ds)\bigg\}.\end{aligned}$$

By (11), $\mathscr{E}_t(\phi)$ is well defined. All such predictable processes will be denoted by $\mathscr{G}$.

PROPOSITION 4.2. *Let $0 < c_0 \leq F \leq c_1$ be a random variable on $\Omega$. Then for some $\phi \in \mathscr{G}$,*

$$\mathbb{E}(F|\mathcal{F}_t) = \mathbb{E}F \cdot \mathscr{E}_t(\phi), \qquad \forall t \in [0,1].$$

*More precisely,*

(16) $$\phi(u,t) = \frac{{}^p D_{(u,t)} F}{\mathbb{E}(F|\mathcal{F}_{t-})}, \qquad \pi \times \mathbb{P}\text{-a.s.}$$

PROOF. By Theorem 3.4, we have

$$M_t := \mathbb{E}(F|\mathcal{F}_t) = \mathbb{E}F + \int_0^t \int_{\mathbb{U}} {}^p D_{(u,s)} F \tilde{\mu}(du,ds).$$

Let $\phi$ be given by (16). Then, it is clear by Theorem 3.4 that $\phi \in \mathbb{L}^2_{\mathcal{P}} \cap \mathbb{L}^\infty_{\mathcal{P}}$ and (14) holds. For (13), it only needs to notice that by (7) and (6)

$$\phi(u,t) = \frac{\mathbb{E}(F \circ \varepsilon^+_{(u,t)}|\mathcal{F}_{t-})}{\mathbb{E}(F|\mathcal{F}_{t-})} - 1$$

$$\geq \frac{c_0}{c_1} - 1, \qquad \pi \times \mathbb{P}\text{-a.s.}$$

On the other hand, if we define

$$X_t := \int_0^t \int_{\mathbb{U}} \phi(u,s)\tilde{\mu}(du,ds),$$

then by $\phi \in \mathbb{L}^2_{\mathcal{P}}$, $X$ is a square-integrable $\mathcal{F}_t^{\mathbb{P}}$-martingale, $\Delta X_s \geq \frac{c_0}{c_1} - 1$, and

$$M_t = \mathbb{E}F + \int_0^t M_{s-} dX_s.$$

By [14], page 84, Theorem 37, we have

$$M_t = \mathbb{E}F \cdot \exp\{X_t\} \cdot \prod_{0<s\leq t}[(1+\Delta X_s)\cdot\exp\{-\Delta X_s\}]$$

(17)
$$= \mathbb{E}F \cdot \exp\left\{X_t + \sum_{0<s\leq t}[\log(1+\Delta X_s) - \Delta X_s]\right\}.$$

Observing that by (11)

$$\sum_{0<s\leq t}[\log(1+\Delta X_s) - \Delta X_s] = \int_0^t\int_{\mathbb{U}}[\log(1+\phi(u,s)) - \phi(u,s)]\mu(du,ds)$$

$$= \int_0^t\int_{\mathbb{U}}[\log(1+\phi(u,s)) - \phi(u,s)]\tilde{\mu}(du,ds)$$

$$+ \int_0^t\int_{\mathbb{U}}[\log(1+\phi(u,s)) - \phi(u,s)]\pi(du,ds),$$

we obtain by substituting this into (17)

$$\mathbb{E}(F|\mathcal{F}_t) = M_t = \mathbb{E}F \cdot \mathscr{E}_t(\phi),$$

which then implies that $t \mapsto \mathscr{E}_t(\phi)$ is a square integrable $\mathcal{F}_t^{\mathbb{P}}$-martingale, and so $\phi \in \mathscr{G}$. □

PROPOSITION 4.3. *For $\phi \in \mathscr{G}$, define a new probability measure on $(\Omega, \mathcal{F}^{\mathbb{P}})$ by*

(18)
$$d\mathbb{P}_\phi := \mathscr{E}_1(\phi)\,d\mathbb{P},$$

*then for any $\psi \in \mathbb{L}_{\mathcal{P}}^2$,*

$$t \mapsto \int_0^t\int_{\mathbb{U}}\psi(u,s)\tilde{\mu}(du,ds) - \int_0^t\int_{\mathbb{U}}\psi(u,s)\phi(u,s)\pi(du,ds)$$

*is a square integrable $\mathcal{F}_t^{\mathbb{P}}$-martingale under $\mathbb{P}_\phi$.*

PROOF. Note that by Itô's formula, $\mathscr{E}_t(\phi)$ solves the following linear equation:

$$\mathscr{E}_t(\phi) = 1 + \int_0^t\int_{\mathbb{U}}\mathscr{E}_{s-}(\phi) \cdot \phi(u,s)\tilde{\mu}(du,ds).$$

If we put $Z_t := \int_0^t\int_{\mathbb{U}}\psi(u,s)\tilde{\mu}(du,ds)$, then

$$\langle Z, \mathscr{E}(\phi)\rangle_t = \int_0^t\int_{\mathbb{U}}\mathscr{E}_{s-}(\phi) \cdot \psi(u,s)\phi(u,s)\pi(du,ds).$$

By the Meyer–Girsanov theorem (cf. [14], page 133, Theorem 36), we know that
$$t \mapsto Z_t - \int_0^t \frac{1}{\mathscr{E}_{s-}(\phi)} \, d\langle Z, \mathscr{E}(\phi) \rangle_s$$
is a square integrable $\mathcal{F}_t^{\mathbb{P}}$-martingale under $\mathbb{P}_\phi$. The result follows. $\square$

We may prove the following representation formula.

THEOREM 4.4. *Let $F$ be a bounded random variable on $\Omega$. Then*
$$\text{(19)} \qquad -\log \mathbb{E}(e^{-F}) = \inf_{\phi \in \mathscr{G}} \mathbb{E}^{\mathbb{P}_\phi}(F + L(\phi)),$$
*where $\mathbb{P}_\phi$ is defined by (18), and*
$$\text{(20)} \quad L(\phi) := \int_0^1 \int_{\mathbb{U}} [(1+\phi(u,s))\log(1+\phi(u,s)) - \phi(u,s)] \pi(du,ds)$$
*is well defined by (12). Moreover, the infimum is uniquely attained at some $\phi \in \mathscr{G}$.*

PROOF. For any $\phi \in \mathscr{G}$, by Jensen's inequality, we have
$$-\log \mathbb{E}(e^{-F}) = -\log \mathbb{E}^{\mathbb{P}_\phi}(e^{-F - \log(d\mathbb{P}_\phi/d\mathbb{P})})$$
$$\leq \mathbb{E}^{\mathbb{P}_\phi}(F) + R(\mathbb{P}_\phi \| \mathbb{P}).$$

By (18) and (15), we have
$$R(\mathbb{P}_\phi \| \mathbb{P}) = \mathbb{E}^{\mathbb{P}_\phi} \bigg( \int_0^1 \int_{\mathbb{U}} \log(1+\phi(u,s)) \tilde{\mu}(du,ds)$$
$$+ \int_0^1 \int_{\mathbb{U}} [\log(1+\phi(u,s)) - \phi(u,s)] \pi(du,ds) \bigg).$$

By Proposition 4.3, we know that
$$t \mapsto \int_0^1 \int_{\mathbb{U}} \log(1+\phi(u,s)) \tilde{\mu}(du,ds) - \int_0^1 \int_{\mathbb{U}} \phi(u,s) \log(1+\phi(u,s)) \pi(du,ds)$$
is a square integrable $\mathcal{F}_t^{\mathbb{P}}$-martingale under $\mathbb{P}_\phi$. Hence, by (12), (14) and Hölder's inequality,
$$R(\mathbb{P}_\phi \| \mathbb{P}) = \mathbb{E}^{\mathbb{P}_\phi} \bigg( \int_0^1 \int_{\mathbb{U}} [(1+\phi(u,s))\log(1+\phi(u,s)) - \phi(u,s)] \pi(du,ds) \bigg)$$
$$\leq \mathbb{E} \bigg( \mathscr{E}_1(\phi) \cdot \int_0^1 \int_{\mathbb{U}} |\phi(u,s)|^2 \pi(du,ds) \bigg) < +\infty.$$

Thus, the upper bound is obtained.

For the lower bound, by Proposition 4.2, there exists a $\phi \in \mathscr{G}$ such that $\mathbb{E}(e^{-F}) = e^{-F} \cdot \mathscr{E}_1^{-1}(\phi)$. Thus, we have

$$-\log \mathbb{E}(e^{-F}) = \mathbb{E}^{\mathbb{P}_\phi}(F) + R(\mathbb{P}_\phi \| \mathbb{P}) = \mathbb{E}^{\mathbb{P}_\phi}(F + L(\phi)).$$

The uniqueness follows from that when the infimum is attained, then Jensen's inequality becomes an equality. The proof is thus complete. $\square$

We now turn to proving another representation like (1) about the formula (19), which was conjectured by Dupuis and Ellis [4], page 122. For further discussions, we need to consider a noncompact locally compact connected complete metric space $(\mathbb{U}, \rho)$, and assume that:

(H1) For each $\phi \in \mathscr{G}$, there exists an invertible transformation with respect to $u$,

$$\gamma_\phi : \mathbb{U} \times [0,1] \times \Omega \to \mathbb{U}, \qquad \mathbb{U} \ni u \mapsto \gamma_\phi(u, t, \omega) \in \mathbb{U},$$

such that
   (i) $\gamma_\phi, \gamma_\phi^{-1} \in \mathcal{B}(\mathbb{U}) \times \mathcal{P}/\mathcal{B}(\mathbb{U})$;
   (ii) $\nu \circ \gamma_\phi^{-1} = (1 + \phi) \cdot \nu$, that is, for $(ds \times d\mathbb{P})$-almost all $(s, \omega) \in [0,1] \times \Omega$ and any bounded measurable function $f$ on $\mathbb{U}$

$$\int_{\mathbb{U}} f(\gamma_\phi(u, s, \omega)) \nu(du) = \int_{\mathbb{U}} f(u) \cdot (1 + \phi(u, s, \omega)) \nu(du);$$

   (iii) for each $t \in [0, 1]$, $\gamma_{\phi_t}|_{[0,t]} = \gamma_\phi|_{[0,t]}$ and $\gamma_{\phi_t}^{-1}|_{[0,t]} = \gamma_\phi^{-1}|_{[0,t]}$, where $\phi_t := \phi \cdot 1_{[0,t]}$.

(H2) Let $\phi, \phi_n \in \mathscr{G}$ satisfy $-1 < c_0 \leq \phi, \phi_n \leq c_1$. If $\phi_n$ converges to $\phi$ in $\mathbb{L}^2_{\mathcal{P}}$, then there is a subsequence $n_k$ (still denoted by $n$) such that for $(\pi \times \mathbb{P})$-almost all $(u, s, \omega)$,

$$\lim_{n \to \infty} \rho(\gamma_{\phi_n}(u, s, \omega), \gamma_\phi(u, s, \omega)) = \lim_{n \to \infty} \rho(\gamma_{\phi_n}^{-1}(u, s, \omega), \gamma_\phi^{-1}(u, s, \omega)) = 0,$$

where $\rho$ is the metric on $\mathbb{U}$.

REMARK 4.5. The invertibility is understood in the measure sense, that is,

$$\gamma_\phi(\gamma_\phi^{-1}(u, t, \omega), t, \omega) = \gamma_\phi^{-1}(\gamma_\phi(u, t, \omega), t, \omega) = u, \qquad \pi \times \mathbb{P}\text{-a.s.}$$

However, by a suitable redefinition procedure, one may assume that the above identities hold for all $(u, t, \omega) \in \mathbb{U} \times [0, 1] \times \Omega$. In fact, for some $(dt \times \mathbb{P})$-null set $A \in \mathcal{P}$ and each $(t, \omega) \notin A$, there exists a $\nu$-null set $N_{(t,\omega)} \in \mathcal{B}(\mathbb{U})$ such that $\gamma_\phi(\cdot, t, \omega)$ is a one-to-one and onto mapping on $N_{(t,\omega)}^c$. Thus, $\gamma_\phi^+ := \gamma_\phi$ and $\gamma_\phi^- := \gamma_\phi^{-1}$ may be redefined as follows:

$$\tilde{\gamma}_\phi^\pm(u, t, \omega) := \begin{cases} \gamma_\phi^\pm(u, t, \omega), & \text{if } (t, \omega) \in A^c \text{ and } u \in N_{(t,\omega)}^c, \\ u, & \text{otherwise.} \end{cases}$$

In the sequel, we still use $\gamma_\phi^\pm$ to denote these redefinitions.

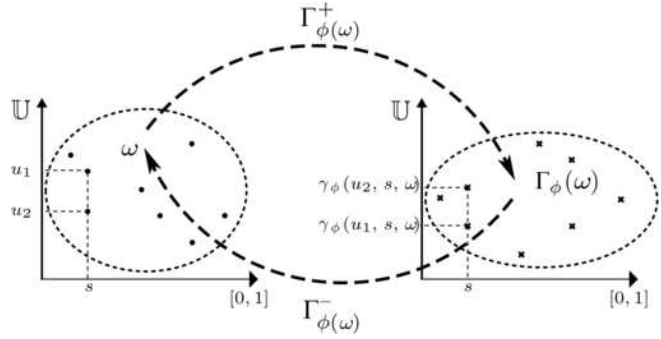

FIG. 1. *Transformation:* $\Gamma_\phi^\pm$.

The above constructed $\gamma_\phi^\pm$ induce predictable transformations $\Gamma_\phi^\pm$ on $\Omega$ as follows (see Figure 1):

$$(21) \qquad \omega \mapsto \Gamma_{\phi(\omega)}^\pm(\omega) := (\gamma_\phi^\pm(\omega))_*(\omega),$$

where $(\gamma_\phi^\pm)_*(\omega)$ denote the image measures of $\omega$ under transformations $(u,t) \mapsto (\gamma_\phi^\pm(u,t,\omega), t)$. In particular, for each $\omega \in \Omega$,

$$(22) \qquad \Gamma_{\phi(\omega)}^+(\Gamma_{\phi(\omega)}^-(\omega)) = \Gamma_{\phi(\omega)}^-(\Gamma_{\phi(\omega)}^+(\omega)) = \omega.$$

In what follows, we sometimes simply write $\Gamma_{\phi(\omega)}^+$ (resp. $\Gamma_{\phi(\omega)}^-$) as $\Gamma_\phi$ (resp. $\Gamma_\phi^-$). The following Girsanov theorem can be found in [5], page 165.

THEOREM 4.6. *Assume* (H1). *For any $\phi \in \mathscr{G}$, the mapping $\omega \mapsto \mu_{\Gamma_\phi(\omega)}$ is still a Poisson random measure under $\mathbb{P}_\phi$ with the same intensity measure $\nu$, where $\mathbb{P}_\phi$ is defined by (18). In particular,*

$$\mathbb{P}_\phi \circ (\Gamma_\phi)^{-1} = \mathbb{P},$$

*where $\mathbb{P}_\phi \circ (\Gamma_\phi)^{-1}$ denotes the image measure or distribution of $\omega \mapsto \Gamma_\phi(\omega)$ under $\mathbb{P}_\phi$.*

We now prepare several lemmas for later use. The following lemma is direct.

LEMMA 4.7. *Assume* (H2). *Then $\mathbb{P} \circ (\Gamma_{\phi_n}^\pm)^{-1}$ weakly converges to $\mathbb{P} \circ (\Gamma_\phi^\pm)^{-1}$ as $n \to \infty$.*

PROOF. It only needs to prove that for $\mathbb{P}$-almost all $\omega \in \Omega$, $\Gamma_{\phi_n}^\pm(\omega)$ converges to $\Gamma_\phi^\pm(\omega)$ with respect to the weak topology defined by (5). That is, for any $f \in C_c(\mathbb{U} \times [0,1])$,

$$\langle f, \mu_{\Gamma^\pm_{\phi_n}(\omega)} \rangle \to \langle f, \mu_{\Gamma^\pm_\phi(\omega)} \rangle.$$

As in Remark 4.5, by (H2), we assume that for $\mathbb{P}$-almost all $\omega$ and all $(u,t) \in \mathbb{U} \times [0,1]$

$$\lim_{n\to\infty} \rho(\gamma^\pm_{\phi_n}(u,t,\omega), \gamma^\pm_\phi(u,t,\omega)) = 0.$$

Since $f$ has compact support in $\mathbb{U} \times [0,1]$, we have

$$\langle f, \mu_{\Gamma^\pm_{\phi_n}(\omega)} \rangle = \sum_{(u,t) \in \mathrm{supp}(\omega)} f(\gamma^\pm_{\phi_n}(u,t,\omega), t)$$

$$\to \sum_{(u,t) \in \mathrm{supp}(\omega)} f(\gamma^\pm_\phi(u,t,\omega), t)$$

$$= \langle f, \mu_{\Gamma^\pm_\phi(\omega)} \rangle.$$

The result follows. $\square$

We introduce the following subclasses of $\mathscr{G}$ and $\mathscr{C}$: For $-1 < c_0 \leq 0$ and $c_1 > 0$, $\phi \in \mathscr{G}^{c_1}_{c_0} \subset \mathscr{G}$ or $\phi \in \mathscr{C}^{c_1}_{c_0} \subset \mathscr{C}$ if

$$c_0 \leq \phi \leq c_1.$$

It is noticed that $\mathscr{C}^{c_1}_{c_0} \subset \mathscr{G}^{c_1}_{c_0}$ by (3).

LEMMA 4.8. *Let $-1 < c_0 \leq 0$ and $c_1 > 0$. For any $\phi \in \mathscr{G}^{c_1}_{c_0}$, there exists a sequence $\phi_n \in \mathscr{C}^{c_1}_{c_0}$ such that*

(23) $$\lim_{n\to\infty} \int_0^1 \int_\mathbb{U} \mathbb{E}|\phi_n(u,t) - \phi(u,t)|^2 \pi(du, dt) = 0,$$

*and*

(24) $$\lim_{n\to\infty} \mathbb{E}|L(\phi_n) - L(\phi)| = 0.$$

PROOF. As in the construction in Lemma 2.2, it is easy to find the desired $\phi_n$. As for the limit (24), it follows from the construction of $\phi_n$ in Lemma 2.2, (12) and the dominated convergence theorem. $\square$

LEMMA 4.9. *Assume* (H1). *Let $g$ be a bounded $\mathcal{F}^\mathbb{P}_t$-measurable function. Then for any $\phi \in \mathscr{G}$,*

$$g(\Gamma^\pm_\phi(\omega)) = g(\Gamma^\pm_{\phi_t}(\omega)), \qquad \mathbb{P}\text{-}a.a.\ \omega,$$

*where $\phi_t = \phi \cdot 1_{[0,t]}$.*

PROOF. By the monotone class theorem, it is enough to consider cylindrical function $g$ with the following form:

$$g(\omega) = h(\langle f_1, \mu_\omega\rangle, \ldots, \langle f_n, \mu_\omega\rangle), \qquad h \in C_c^\infty(\mathbb{R}^n), f_i \in C_c(\mathbb{U} \times [0,t]).$$

For this type $g$, the desired equality follows by direct calculations and (iii) of (H1). □

The following lemma is crucial for the proof of Theorem 4.11 below. The main idea comes from [1, 3] (see also [21]).

LEMMA 4.10. *Assume* (H1). *Let* $-1 < c_0 \leq 0$ *and* $c_1 > 0$. *For any* $\phi \in \mathscr{C}_{c_0}^{c_1}$, *there are two* $\tilde{\phi}, \hat{\phi} \in \mathscr{C}_{c_0}^{c_1}$ *such that for any bounded random variable* $F$ *on* $\Omega$

(25) $$\mathbb{E}^{\mathbb{P}\tilde{\phi}}(F + L(\tilde{\phi})) = \mathbb{E}(F \circ \Gamma_\phi^- + L(\phi)),$$

(26) $$\mathbb{E}^{\mathbb{P}\phi}(F + L(\phi)) = \mathbb{E}(F \circ \Gamma_{\hat{\phi}}^- + L(\hat{\phi})),$$

*where the functional* $L$ *is defined by (20). Moreover,*

(27) $$R(\mathbb{P} \circ (\Gamma_\phi^-)^{-1} \| \mathbb{P}) = \mathbb{E}^{\mathbb{P}\tilde{\phi}}(L(\tilde{\phi})) = \mathbb{E}L(\phi).$$

PROOF. Let $\phi \in \mathscr{C}_{c_0}^{c_1}$ have the form

$$\phi(u,t,\omega) := \sum_{i=0}^{n} 1_{(t_i, t_{i+1}]}(t) \cdot g_i(u,\omega), \qquad g_i \in \mathcal{B}(\mathbb{U}) \times \mathcal{F}_{t_i}.$$

Let us construct $\tilde{g}_i$ as follows:

$$\tilde{g}_0(u,\omega) = g_0(u,\omega)$$

and for $i = 1, 2, \ldots, n-1$,

$$\tilde{g}_i(u,\omega) = g_i(u, \Gamma_{\tilde{\phi}_i(\omega)}(\omega)),$$

where

$$\tilde{\phi}_i(u,t,\omega) := \sum_{j=0}^{i-1} 1_{(t_j, t_{j+1}]}(t) \cdot \tilde{g}_j(u,\omega).$$

Finally, we let

$$\tilde{\phi}(u,t,\omega) := \tilde{\phi}_n(u,t,\omega).$$

From the construction, it is clear that $\tilde{\phi} \in \mathscr{C}_{c_0}^{c_1}$. Moreover, it is not hard to verify by Lemma 4.9 and induction that $\tilde{\phi}$ satisfies

(28) $$\tilde{\phi}(u,t,\omega) = \phi(u,t,\Gamma_{\tilde{\phi}(\omega)}(\omega)).$$

Similarly, one may construct $\hat{\phi} \in \mathscr{C}_{c_0}^{c_1}$ such that

$$\hat{\phi}(u,t,\omega) := \phi(u,t,\Gamma^-_{\hat{\phi}(\omega)}(\omega)).$$

As above, by induction, Lemma 4.9 and (22), one can verify

(29) $$\phi(u,t,\omega) = \hat{\phi}(u,t,\Gamma_{\phi(\omega)}(\omega)).$$

Now by Theorem 4.6, we have

(30) $$\mathbb{P}_{\tilde{\phi}} \circ (\Gamma_{\tilde{\phi}})^{-1} = \mathbb{P} = \mathbb{P}_\phi \circ (\Gamma_\phi)^{-1}.$$

Hence, we obtain by (22) and (28)

$$\mathbb{E}^{\mathbb{P}_{\tilde{\phi}}}(F + L(\tilde{\phi})) = \mathbb{E}^{\mathbb{P}_{\tilde{\phi}}}(F(\Gamma^-_{\phi(\Gamma_{\tilde{\phi}})}(\Gamma_{\tilde{\phi}}(\cdot))) + L(\phi(\Gamma_{\tilde{\phi}})))$$

$$= \mathbb{E}(F \circ \Gamma^-_\phi + L(\phi)),$$

as well as by (22) and (29)

$$\mathbb{E}^{\mathbb{P}_\phi}(F + L(\phi)) = \mathbb{E}^{\mathbb{P}_\phi}(F(\Gamma^-_{\hat{\phi}(\Gamma_\phi)}(\Gamma_\phi(\cdot))) + L(\hat{\phi}(\Gamma_\phi)))$$

$$= \mathbb{E}(F \circ \Gamma^-_{\hat{\phi}} + L(\hat{\phi})).$$

Moreover, by (30), (28) and (22), we also have

$$\mathbb{P} \circ (\Gamma^-_\phi)^{-1} = \mathbb{P}_{\tilde{\phi}}$$

and so,

$$R(\mathbb{P} \circ (\Gamma^-_\phi)^{-1} \| \mathbb{P}) = R(\mathbb{P}_{\tilde{\phi}} \| \mathbb{P}) = \mathbb{E}^{\mathbb{P}_{\tilde{\phi}}}(L(\tilde{\phi})) = \mathbb{E}^{\mathbb{P}_{\tilde{\phi}}}(L(\phi(\Gamma_{\tilde{\phi}}))) = \mathbb{E}L(\phi).$$

The proof is complete. □

We are now in a position to prove our main result in the present paper.

THEOREM 4.11. *Assume that* (H1) *and* (H2) *hold. Let $F$ be any bounded random variable on $\Omega$. Then*

$$-\log \mathbb{E}(e^{-F}) = \inf_{\phi \in \mathscr{G}} \mathbb{E}(F \circ \Gamma^-_\phi + L(\phi))$$

$$= \inf_{\phi \in \mathscr{C}_\alpha^\beta} \mathbb{E}(F \circ \Gamma^-_\phi + L(\phi)),$$

*where $L(\phi)$ and $\Gamma^-_\phi$ are defined by (20) and (21) respectively, and*

(31) $$\alpha := e^{-2\|F\|_\infty} - 1, \qquad \beta := 1 + e^{2\|F\|_\infty}.$$

PROOF. (*Upper bound*) Let $\phi \in \mathscr{G}$. Then $\phi \in \mathscr{G}_{c_0}^{c_1}$ for some $c_0 \in (-1, 0]$, $c_1 > 0$. Let $\phi_n \in \mathscr{C}_{c_0}^{c_1}$ be as in Lemma 4.8. Let $\tilde{\phi}_n \in \mathscr{C}_{c_0}^{c_1}$ be the corresponding one constructed in Lemma 4.10. Then, by (19) and (25),

$$(32) \quad -\log \mathbb{E}(e^{-F}) \leq \mathbb{E}^{\mathbb{P}\tilde{\phi}_n}(F + L(\tilde{\phi}_n)) = \mathbb{E}(F \circ \Gamma_{\phi_n}^- + L(\phi_n)).$$

Noting that by (27) and (24)

$$\sup_n R(\mathbb{P} \circ (\Gamma_{\phi_n}^-)^{-1} \| \mathbb{P}) = \sup_n \mathbb{E}L(\phi_n) < +\infty,$$

we have by Lemma 4.7 and (ii) of Lemma 2.5

$$\lim_{n \to \infty} \mathbb{E}(F \circ \Gamma_{\phi_n}^-) = \mathbb{E}(F \circ \Gamma_\phi^-).$$

Hence, by (32) and (24),

$$-\log \mathbb{E}(e^{-F}) \leq \mathbb{E}(F \circ \Gamma_\phi + L(\phi)),$$

which gives the upper bound.

Moreover, by the lower semi-continuity of $R(\cdot \| \mathbb{P})$ (cf. [4], Lemma 1.4.3), we also have

$$(33) \quad \begin{aligned} R(\mathbb{P} \circ (\Gamma_\phi^-)^{-1} \| \mathbb{P}) &\leq \varliminf_{n \to \infty} R(\mathbb{P} \circ (\Gamma_{\phi_n}^-)^{-1} \| \mathbb{P}) \\ &= \lim_{n \to \infty} \mathbb{E}L(\phi_n) \leq \mathbb{E}L(\phi), \qquad \text{for all } \phi \in \mathscr{G}. \end{aligned}$$

(*Lower bound*) We divide the proof into two steps.

(Step 1): First of all, let $F \in \mathcal{C}$ have the following form:

$$F(\omega) = g(\langle f_1, \mu_\omega \rangle, \ldots, \langle f_n, \mu_\omega \rangle), \qquad g \in C_c^\infty(\mathbb{R}^n), f_i \in C_c(\mathbb{U} \times [0,1]).$$

Then, by (6) and a simple calculation, we have

$$(34) \quad |D_{(u,t)} e^{-F(\omega)}| = |e^{-F(\varepsilon_{(u,t)}^+ \omega)} - e^{-F(\omega)}| \leq C \sum_{i=1}^n |f_i(u,t)|,$$

where $C$ is independent of $(u, t, \omega)$.

Set

$$\phi(u, t) := \frac{{}^p D_{(u,t)} e^{-F}}{\mathbb{E}(e^{-F}|\mathcal{F}_{t-})}.$$

It is clear by Proposition 4.2 that $\phi \in \mathscr{G}_\alpha^\beta$, where $\alpha, \beta$ are given by (31). Let $\phi_n \in \mathscr{C}_\alpha^\beta$ be as in Lemma 4.8. By (34) and the construction of $\phi_n$, there exists a $U \subset \mathbb{U}$ with $\nu(U) < +\infty$ such that for all $n \in \mathbb{N}$

$$(35) \quad |\phi_n(u, t, \omega)| \leq C \cdot 1_U(u), \qquad \pi \times \mathbb{P}\text{-a.e.}$$

By limits (23), (24) and extracting a subsequence if necessary, we may further assume that
$$\phi_n \to \phi, \qquad \pi \times \mathbb{P}\text{-a.e.}$$
and
$$L(\phi_n) \to L(\phi), \qquad \tilde{\mu}(\phi_n) \to \tilde{\mu}(\phi), \qquad \mathbb{P}\text{-a.e.}$$

By (35) and the dominated convergence theorem, we have

(36) $$\mathbb{E}^{\mathbb{P}_{\phi_n}}(F + L(\phi_n)) \to \mathbb{E}^{\mathbb{P}_\phi}(F + L(\phi)) \qquad \text{as } n \to \infty.$$

Moreover, by Proposition 4.2, we have
$$e^{-F} = \mathbb{E}(e^{-F})\mathscr{E}_1(\phi).$$

So, by (36) and (26), we have for any $\varepsilon > 0$ and $n$ large enough
$$-\log \mathbb{E}(e^{-F}) = \mathbb{E}^{\mathbb{P}_\phi}(F) + R(\mathbb{P}_\phi \| \mathbb{P})$$
$$= \mathbb{E}^{\mathbb{P}_\phi}(F + L(\phi))$$
$$\geq \mathbb{E}^{\mathbb{P}_{\phi_n}}(F + L(\phi_n)) - \varepsilon$$
$$= \mathbb{E}(F \circ \Gamma^-_{\hat{\phi}_n} + L(\hat{\phi}_n)) - \varepsilon.$$

The lower bound now follows by $\hat{\phi}_n \in \mathscr{C}^\beta_\alpha$ (see Lemma 4.10).

(Step 2): For any bounded random variable $F$ on $(\Omega, \mathcal{F})$, by Lemma 2.6, there exists a sequence $F_n \in \mathcal{C}$ such that

(37) $$\sup_n \|F_n\|_\infty \leq \|F\|_\infty$$

and
$$\lim_{n \to \infty} F_n = F, \qquad \mathbb{P}\text{-a.s.}$$

For any $\varepsilon > 0$ and $F_n$, by Step 1 and (37), there exists a $\phi_n \in \mathscr{C}^\beta_\alpha$, where $\alpha, \beta$ are given by (31), such that

(38) $$-\log \mathbb{E}(e^{-F_n}) \geq \mathbb{E}(F_n \circ \Gamma^-_{\phi_n} + L(\phi_n)) - \varepsilon.$$

In view of (33), (38) and (37), we have
$$\sup_n R(\mathbb{P} \circ (\Gamma^-_{\phi_n})^{-1} \| \mathbb{P}) \leq \sup_n \mathbb{E}L(\phi_n) < +\infty.$$

Therefore, by (i) of Lemma 2.5,
$$\lim_{n \to \infty} \mathbb{E}|F_n \circ \Gamma^-_{\phi_n} - F \circ \Gamma^-_{\phi_n}| = 0.$$

Using the dominated convergence theorem to the left-hand side of (38) gives that, for sufficiently large $n$,

$$-\log \mathbb{E}(e^{-F}) \geq \mathbb{E}(F \circ \Gamma_{\phi_n}^- + L(\phi_n)) - 2\varepsilon.$$

Since $\phi_n \in \mathscr{C}_\alpha^\beta$ and $\varepsilon$ is arbitrary, we thus complete the proof of the lower bound. □

REMARK 4.12. By the same argument as in the proof of [1], Theorem 5.1, the $F$ in Theorem 4.11 can be any random variable bounded from above.

**5. (H1)–(H2) and mass transportation problem.** In this section we give a more concrete description for (H1)–(H2). Let $(\mathbb{U}, \rho)$ be a locally compact complete metric space, and $\nu$ a $\sigma$-finite and infinite measure on $(\mathbb{U}, \mathcal{B}(\mathbb{U}))$. Let $\mathscr{U}$ be the set of all positive measurable functions on $\mathbb{U}$ bounded from above and also from below.

QUESTION. Under what constraints, for each $\phi \in \mathscr{U}$, does there exist a unique invertible measurable transformation $\gamma_\phi$ on $\mathbb{U}$ such that $\nu \circ \gamma_\phi^{-1} = \phi \cdot \nu$, that is,

$$(39) \qquad \int_\mathbb{U} f(\gamma_\phi(u))\nu(du) = \int_\mathbb{U} f(u)\phi(u)\nu(du), \qquad \forall f \in C_c(\mathbb{U})?$$

Moreover, for $0 < C_0 \leq \phi, \phi_n \leq C_1$, if $\phi_n$ converges $\nu$-a.e. to $\phi$, does it hold that

$$(40) \qquad \lim_{n \to \infty} \rho(\gamma_{\phi_n}(u), \gamma_\phi(u)) = \lim_{n \to \infty} \rho(\gamma_{\phi_n}^{-1}(u), \gamma_\phi^{-1}(u)) = 0, \qquad \nu\text{-a.a. } u?$$

Obviously, if this question has a solution, then (H1) and (H2) are satisfied. We remark that the required predictability follows from continuous dependence (40) with respect to $\phi$. In the classical problem of optimal mass transportation, the constraint is given by minimizing the following cost functional (cf. [20]):

$$\inf_{\gamma_\phi} \int_\mathbb{U} c(\rho(u, \gamma_\phi(u)))\nu(du),$$

where $c$ is a convex function on $\mathbb{R}_+$.

Let us look at the case of $\mathbb{U} = \mathbb{R}^d$ and $\nu(dx) = \theta(x)\,dx$. It is clear that (39) can be reduced to

$$\theta(\gamma_\phi^{-1}(x)) \det(\nabla \gamma_\phi^{-1}(x)) = \phi(x)\theta(x),$$

where $\nabla$ denotes the gradient. If we further require that $\gamma_\phi^{-1}(x)$ is the gradient of some strictly convex function $h(x)$, then we need to solve the following classical Monge–Ampère equation:

$$\theta(\nabla h(x)) \det(\nabla^2 h(x)) = \phi(x)\theta(x).$$

For this equation, there are many literatures to study it, for example, see [20] and references therein. Since our problem is looser, we can find an easy solution when $\mathbb{U} = \mathbb{R}^d$.

Let us first see the one-dimensional case. Let $\nu$ have full support and no atoms. There are three possibilities:

(1) $\nu([0, +\infty)) = \nu((-\infty, 0]) = +\infty$,
(2) $\nu([0, +\infty)) = +\infty$, $\nu((-\infty, 0]) < +\infty$,
(3) $\nu([0, +\infty)) < +\infty$, $\nu((-\infty, 0]) = +\infty$.

It suffices to consider the first case. The others are analogous. Let $\phi \in \mathscr{U}$. In the first case, note that for $x \geq 0$

$$\Phi_+(x) := \int_0^x \phi(u)\nu(du) \quad \text{and} \quad \Phi_-(x) := \int_{-x}^0 \phi(u)\nu(du)$$

are strictly increasing continuous functions on $[0, +\infty)$, and $\Phi_\pm(+\infty) = +\infty$.

Define for $x \geq 0$

$$\gamma_\phi(x) := \Phi_+^{-1}(\nu([0, x])), \qquad \gamma_\phi(-x) := -\Phi_-^{-1}(\nu([-x, 0])).$$

It is clear that $\gamma_\phi$ is an invertible continuous transformation of $\mathbb{R}$, and for any $a < b$,

(41) $$\nu([a, b]) = \int_{\gamma_\phi(a)}^{\gamma_\phi(b)} \phi(u)\nu(du),$$

which means $\nu \circ \gamma_\phi^{-1} = \phi \cdot \nu$.

We now verify (40). For $0 < C_0 \leq \phi, \phi_n \leq C_1$, assume that $\phi_n$ converges $\nu$-a.s. to $\phi$. Noticing that by (41)

$$\int_0^{\gamma_{\phi_n}(x)} \phi_n(u)\nu(du) = \int_0^{\gamma_\phi(x)} \phi(u)\nu(du),$$

we have by the dominated convergence theorem

$$\left| \int_{\gamma_\phi(x)}^{\gamma_{\phi_n}(x)} \nu(du) \right| \leq \frac{1}{C_0} \left| \int_{\gamma_\phi(x)}^{\gamma_{\phi_n}(x)} \phi_n(u)\nu(du) \right|$$

$$= \frac{1}{C_0} \left| \int_0^{\gamma_\phi(x)} (\phi_n(u) - \phi(u))\nu(du) \right| \to 0.$$

Since $\nu$ has full support in $\mathbb{R}$, it follows that

$$\lim_{n \to \infty} |\gamma_{\phi_n}(x) - \gamma_\phi(x)| = 0.$$

Similarly,

$$\lim_{n \to \infty} |\gamma_{\phi_n}^{-1}(x) - \gamma_\phi^{-1}(x)| = \lim_{n \to \infty} \left| \int_0^{\gamma_\phi^{-1}(x)} (\phi_n(u) - \phi(u))\nu(du) \right| = 0.$$

For the multi-dimensional case, we assume that $\nu$ has full support and no charges on $d-1$-dimensional subspaces, and one of the following conditions holds: For any $x_i \in \mathbb{R} - \{0\}, i = 2, \ldots, d$, let $x_i^+ = x_i \vee 0$ and $x_i^- = x_i \wedge 0$,

(1') $\int_0^\infty \int_{x_2^-}^{x_2^+} \cdots \int_{x_d^-}^{x_d^+} d\nu = \infty$, $\int_{-\infty}^0 \int_{x_2^-}^{x_2^+} \cdots \int_{x_d^-}^{x_d^+} d\nu < \infty$;

(2') $\int_0^\infty \int_{x_2^-}^{x_2^+} \cdots \int_{x_d^-}^{x_d^+} d\nu < \infty$, $\int_{-\infty}^0 \int_{x_2^-}^{x_2^+} \cdots \int_{x_d^-}^{x_d^+} d\nu = \infty$;

(3') $\int_0^\infty \int_{x_2^-}^{x_2^+} \cdots \int_{x_d^-}^{x_d^+} d\nu = \infty$, $\int_{-\infty}^0 \int_{x_2^-}^{x_2^+} \cdots \int_{x_d^-}^{x_d^+} d\nu = \infty$.

REMARK 5.1. Let $\theta \geq c_0 > 0$ be a continuous function on $\mathbb{R}^d$. If $\nu(dx) = \theta(x)\,dx$, then (3') holds.

We consider the first case. The others are analogous. Without loss of generality, we assume $d = 2$ and fix a $\phi \in \mathscr{U}$. For $x_1, x_2 \in \mathbb{R}$, with $x_2 \neq 0$ let $\alpha_\phi(x_1, x_2)$ and $\beta_\phi(x_1, x_2)$ be the unique elements in $\mathbb{R}$ such that

$$\int_{-\infty}^{x_1} \int_{x_2^-}^{x_2^+} d\nu = \int_{-\infty}^{\alpha_\phi(x_1,x_2)} \int_{x_2^-}^{x_2^+} \phi\,d\nu$$

and

$$\int_{-\infty}^{\beta_\phi(x_1,x_2)} \int_{x_2^-}^{x_2^+} d\nu = \int_{-\infty}^{x_1} \int_{x_2^-}^{x_2^+} \phi\,d\nu.$$

For $x_2 = 0$, set $\alpha_\phi(x_1, 0) = x_1 = \beta_\phi(x_1, 0)$. By the assumption (1'), $\alpha_\phi$ and $\beta_\phi$ are well defined functions on $\mathbb{R} \times \mathbb{R}$, and $\alpha_\phi(\infty, x_2) = \beta(\infty, x_2) = \infty$.

Thus, we may define for $(x_1, x_2) \in \mathbb{R}^2$

$$\gamma_\phi(x_1, x_2) = (\alpha_\phi(x_1, x_2), x_2)$$

and

$$\gamma_\phi^{-1}(x_1, x_2) = (\beta_\phi(x_1, x_2), x_2).$$

It is easy to see that

$$\alpha_\phi(\beta_\phi(x_1, x_2), x_2) = \beta_\phi(\alpha_\phi(x_1, x_2), x_2) = (x_1, x_2)$$

and

$$\gamma_\phi \circ \gamma_\phi^{-1}(x_1, x_2) = \gamma_\phi^{-1} \circ \gamma_\phi(x_1, x_2) = (x_1, x_2).$$

Let $0 < C_0 \leq \phi, \phi_n \leq C_1$ and $\phi_n$ converge $\nu$-a.s. to $\phi$. As in the one-dimensional case, one can prove

$$\lim_{n \to \infty} |\alpha_{\phi_n}(x_1, x_2) - \alpha_\phi(x_1, x_2)| = 0$$

and

$$\lim_{n\to\infty} |\beta_{\phi_n}(x_1, x_2) - \beta_\phi(x_1, x_2)| = 0.$$

Hence,

$$\lim_{n\to\infty} |\gamma_{\phi_n}(x_1, x_2) - \gamma_\phi(x_1, x_2)| = 0$$

and

$$\lim_{n\to\infty} |\gamma_{\phi_n}^{-1}(x_1, x_2) - \gamma_\phi^{-1}(x_1, x_2)| = 0.$$

Summarizing the above discussions, we obtain the following result by Theorem 4.11 and Remark 4.12 when $\mathbb{U} = \mathbb{R}^d$.

THEOREM 5.2. *Let $\nu$ be a $\sigma$-finite and infinite measure on $\mathbb{R}^d$ with full support in $\mathbb{R}^d$ and without charges on any $d-1$-dimensional subspaces. Assume that one of $(1')$, $(2')$ and $(3')$ holds. Then, for any random variable $F$ on $\Omega$ bounded from above,*

$$-\log \mathbb{E}(e^{-F}) = \inf_{\phi \in \mathscr{G}} \mathbb{E}(F \circ \Gamma_\phi^- + L(\phi)),$$

*where $L(\phi)$ and $\Gamma_\phi^-$ are defined respectively by (20) and (21).*

**Acknowledgments.** The author would like to thank Professor Benjamin Goldys for providing him an excellent environment to work in the University of New South Wales. He is also very grateful to Professor Jiagang Ren for his valuable suggestions.

SCHOOL OF MATHEMATICS AND STATISTICS
UNIVERSITY OF NEW SOUTH WALES
SYDNEY, 2052
AUSTRALIA
AND
DEPARTMENT OF MATHEMATICS
HUAZHONG UNIVERSITY OF SCIENCE AND TECHNOLOGY
WUHAN, HUBEI 430074
P. R. CHINA
E-MAIL: XichengZhang@gmail.com